\documentclass{article}

\usepackage{multicol}
\usepackage{amsmath,amsfonts,amsthm,amssymb}
\usepackage{graphicx}
\usepackage{tikz}
\usepackage{enumitem}
\usepackage{verbatim}
\usepackage{color}
\usepackage{xcolor}
\usepackage[normalem]{ulem}
\usepackage{comment,nicefrac}

\newtheorem{thm}{Theorem}

\newtheorem{cor}[thm]{Corollary}

\newtheorem{qst}[thm]{Question}
\newtheorem{prob}[thm]{Problem}

\def\G{{\Gamma}}

\def\e{{\epsilon}}

\def\p{{\pi}}

\def\t{{\tau}}
\def\Th{{\Theta}}
\def\w{{\omega}}
\def\W{{\Omega}}

\def\cD{{\cal D}}

\def\cG{{\cal G}}

\def\cK{{\cal K}}

\def\cP{{\cal P}}
\def\cQ{{\cal Q}}

\def\cS{{\cal S}}

\def\cZ{{\cal Z}}

\def\sD{{\sf D}}
\def\sH{{\sf H}}
\def\sP{{\sf P}}
\def\sw{{\sf w}}

\def\pr{{\rm Pr}}
\def\Pr{{\bf P}}

\def\zE{{\mathbb E}}
\def\zN{{\mathbb N}}

\def\zZ{{\mathbb Z}}

\def\dav{{\sf dav}}

\def\gcd{{\sf gcd}}
\def\lcm{{\sf lcm}}

\def\rank{{\sf rank}}

\def\sP2P{{$\Pi_2^{\sP}$}}
\def\3sat{{\sf 3{\small SAT}}}
\def\liloh{{\operatorname{o}}}

\def\Phalf{{\Pr_{\nicefrac{1}{2}}}}

\def\rar{{\rightarrow}}

\def\sse{{\subseteq}}

\definecolor{brwn}{RGB}{140, 70, 20}
\definecolor{gren}{RGB}{  0,140, 10}

\definecolor{primo}{RGB}{128,0,128}
\definecolor{segundo}{RGB}{255,92,0}
\definecolor{primox}{RGB}{128,0,128}
\definecolor{segundox}{RGB}{255,92,0}

\definecolor{fir}{RGB}{0,104,180}
\definecolor{sec}{RGB}{200, 0, 50}
\definecolor{thi}{RGB}{70, 100, 30}
\definecolor{fou}{RGB}{140, 70, 20}
\definecolor{fif}{RGB}{100, 0, 250}
\definecolor{darkblue}{RGB}{0,0,250}

\pgfmathsetmacro\Radi{6}
\pgfmathsetmacro\radi{3}

\begin{document}

\title{Thresholds for Zero-Sums with Small Cross Numbers\\
in Abelian Groups} 

\author{
Neal Bushaw 
    \thanks{
        Department of Mathematics and Applied Mathematics, Virginia Commonwealth University, Richmond, Virginia, USA}
    \thanks{
        \texttt{ghurlbert@vcu.edu}}
\and
Glenn Hurlbert 
    \footnotemark[1]
    \thanks{
        \texttt{nobushaw@vcu.edu}}
}

\date{}

\maketitle

%%%%%%%%%%%%%%%%%%%%%%%%%%%%%%%%%%%%%%%%%%%%%%%
%       ABSTRACT
%%%%%%%%%%%%%%%%%%%%%%%%%%%%%%%%%%%%%%%%%%%%%%%
\begin{abstract}
For an additive group $\G$ the sequence $S = (g_1, \ldots, g_t)$ of elements of $\G$ is a {\it zero-sum} sequence if $g_1 + \cdots + g_t = 0_\G$.
The {\it cross number} of $S$ is defined to be the sum $\sum_{i=1}^k 1/|g_i|$, where $|g_i|$ denotes the order of $g_i$ in $\G$.
Call $S$ {\it good} if it contains a zero-sum subsequence with cross number at most 1.
In 1993, Geroldinger proved that if $\G$ is abelian then every length $|\G|$ sequence of its elements is good, generalizing a 1989 result of Lemke and Kleitman that had proved an earlier conjecture of Erd\H{o}s and Lemke.
In 1989 Chung re-proved the Lemke and Kleitman result by applying a theorem of graph pebbling, and in 2005, Elledge and Hurlbert used graph pebbling to re-prove and generalize Geroldinger's result.
Here we use probabilistic theorems from graph pebbling to derive a  threshold version of Geroldinger's theorem for abelian groups of a certain form.
Specifically, we prove that if $p_1, \ldots, p_d$ are (not necessarily distinct) primes and $\G_k$ has the form $\prod_{i=1}^d \zZ_{p_i^k}$ then there is a function $\t=\t(k)$ (which we specify in Theorem \ref{t:BushHurl}) with the following property: if $t-\t\rar\infty$ as $k\rar\infty$ then the probability that $S$ is good in $\G_k$ tends to $1$.
\end{abstract}

\medskip


%%%%%%%%%%%%%%%%%%%%%%%%%%%%%%%%%%%%%%%%%%%%%%%
%       INTRO
%%%%%%%%%%%%%%%%%%%%%%%%%%%%%%%%%%%%%%%%%%%%%%%
\section{Introduction}
\label{s:Intro}

For an additive group $\G$ the sequence $S = (g_1, \ldots, g_t)$ of elements of $\G$ is a {\it zero-sum} sequence if $g_1 + \cdots + g_t = 0_\G$.
The study of zero sums has a long and rich history (see \cite{Chapman}).
The following theorem is considered to be one of the jewels and starting points. 

\begin{thm}[Erd\H os-Ginzburg-Ziv \cite{ErdGinZiv}]
\label{t:EGZ}
For any positive integer $n$, every sequence of $2n-1$ of elements from $\zZ_n$ contains a zero-sum subsequence of length exactly $n$.
\end{thm}

Variations on this theme have arisen over time, including the replacement of the length condition in the conclusion by other properties of interest, such as having small cross number, which is important in the area of Krull monoid factorization (see \cite{Chapman}).
For a sequence $S = (g_1, \ldots, g_t)$ of elements of a finite group $\G$, define the {\it cross number} of $S$ to be the sum $\sum_{i=1}^k 1/|g_i|$, where $|g_i|$ denotes the order of $g_i$ in $\G$.
We call $S$ {\it good} if it contains a zero-sum subsequence with cross number at most one.

In 1987 Erd\H{o}s and Lemke \cite{LemkKlei} conjectured that every sequence $(a_1, \ldots, a_n)$ of $n$ elements of $\zZ_n$ contains a zero-sum subsequence whose sum is at most $\lcm(n, a_1, \ldots, a_n)$.
This conjecture was proven by the following theorem.

\begin{thm}[Lemke-Kleitman \cite{LemkKlei}, Chung \cite{Chung}]
\label{t:LK}
For any positive integer $n$, every sequence of $n$ elements of $\zZ_n$ is good.
\end{thm}

The proof by Lemke and Kleitman used what might be considered as traditional techniques.
However, Chung's proof used an ingenious idea of Lagaria and Saks \cite{Chung} that converted the problem to one of pebbling in graphs, which we will discuss in Section \ref{s:Pebbling}, below.
Geroldinger then generalized Theorem \ref{t:LK} to finite abelian groups, along the lines of Lemke and Kleitman's proof, while Elledge and Hurlbert were able to use the pebbling model to attain the same result.

\begin{thm}[Geroldinger \cite{Geroldinger}, Elledge-Hurlbert \cite{ElleHurl}]
\label{t:Geroldinger}
For any finite abelian group $\G$, every sequence of $|\G|$ elements of $\G$ is good.
\end{thm}

Lemke and Kleitman conjectured in \cite{LemkKlei} that the conclusion of Theorem \ref{t:Geroldinger} holds for nonabelian $\G$ also, although no work seems to have been done on this.

% =================================
\subsection{Our results}
\label{ss:Results}

In this paper we move from guaranteeing the result with absolute certainty to achieving the result with high probability, using the pebbling theorems of Section \ref{s:ThreshPebb} to yield sharp thresholds for the existence of zero sums with small cross number.
In this section, we give a short non-technical introduction to the ideas present, so that we can state our results.
We will then revisit these topics in detail in later sections, giving formal definitions to the undefined, and explaining the subtleties that arise.

We describe a randomized version of our sequences, explain how to use pebbling to study  zero-sum sequences to pebbling in Section \ref{ss:Pebb0Sums}, and transfer these more traditional tools to the randomized world via the random pebbling model and theorems in Section \ref{s:ThreshPebb}.

Consider a finite abelian group $\G$, and pick a sequence from $\G$ at random -- that is, we select a sequence $(a_1,\ldots,a_t)$ with equal probability among all sequences of the same length; we denote this probability space by $\cS(\G,t)$.
How large must $t$ be to make it more-likely-than-not that our randomly selected sequence contains a zero-sum subsequence?
How large should $t$ to make it likely that our randomly sequence is good?

With this goal in mind, we select our sequence of length $t$ uniformly from $\G$ and consider the event $\zE_{\G}(t)$ containing all such good sequences. Formally, we will examine the function 
$$\Phalf(\G):=\min\left\{t:\pr[\zE_{\G}(t)]\ge\frac{1}{2}\right\}.$$

We will use Bachmann-Landau notation to describe the asymptotics of functions.
%We will frequently give asymptotics of functions using Bachmann-Landau notation, and it is a standard excercise to determine which of these thresholds describe the asymptotics closely enough to be sharp; in the text below, we will emphasize those results which are sharp.
In this paper we prove the following theorem.

\begin{thm}
\label{t:BushHurl}
Let $\G\cong \prod_{i=1}^d \zZ_{p_i^k}$ be a finite abelian group, where $p_1,\ldots,p_d$ are (not necessarily distinct) primes and $k\ge 1$. Then \[\Phalf(\G)\leq k^d\exp\left[\left(\frac{(d+1)!\prod_i\lg p_i}{2}\lg k\right)^{\frac{1}{d+1}} - \left(\frac{d}{d+1}\right)\lg\lg k + O(1)\right].\]
\end{thm}

Given functions $f$ and $g$, we will write that $f\ll g$, or equivalently $f\in \liloh(g)$ whenever $\lim_{n\to\infty}\frac{f(n)}{g(n)}=0$.
A function $t(n)$ is a {\it{threshold for event $E_n$}} if in the sequence of probability spaces $\cD(G_n,f(n))$ the probability of event $E_n$ tends to $1$ whenever $f(n)\gg t(n)$, and the probability of event $E_n$ tends to $0$ whenever $f\ll t(n)$.
We further call $t(n)$ a {\it{sharp threshold}} if whenever $f(n)\ge(1+\varepsilon)t(n)$ the probability of $E_n$ tends to $1$, and whenever $f(n)\le(1-\varepsilon)t(n)$ the probability of $E_n$ tends to $0$.
Often, for $\cG=(G_i)_{i\in\zN}$, one writes $\tau(\cG)$ for the set of threshold functions for the relevant sequence of events $E_n$.
We will refer to any threshold function that has not been shown to be sharp as a {\it{weak threshold}} (this is sometimes called a {\it{coarse threshold}} in the literature).

These thresholds have the origins in the physical phase transitions that we see in nature, where the physical properties of a state of matter transform dramatically when crossing some critical temperature threshold.
There is a long history and tradition of determining thresholds for random events in both pebbling (for solvability and other events) and random graphs (for connectivity, Hamiltonicity, the existence of a giant component, and many other events).
By determining the asymptotics of $\Phalf(\G)$, we are equivalently identifying the order of magnitude of the functions in $\tau(\G)$.

We note that in the case that our groups are a sequence cyclic groups of increasing powers of a fixed prime $p$, Theorem \ref{t:BushHurl} is improved via the stronger threshold result for pebbling paths due to Bushaw and Kettle (Theorem \ref{t:BushKettPath}); we record this simplest case of Theorem \ref{t:BushHurl} below.

\begin{cor}
\label{t:Examplepk}
Let $\G\cong \prod_{i=1}^d \zZ_{p_i}$ be a finite abelian group, where $p_1,\ldots,p_d$ are (not necessarily distinct) primes. Then \[\Phalf(\G)= \Omega\left(ke^{\sqrt{\lg p\lg k}-(\lg\lg k)/2+o(1)}\right).\]
\end{cor}

%%%%%%%%%%%%%%%%%%%%%%%%%%%%%%%%%%%%%%%%%%%%%%%
%       LS PEBBLING
%%%%%%%%%%%%%%%%%%%%%%%%%%%%%%%%%%%%%%%%%%%%%%%
\section{Lagarias-Saks Pebbling} 
\label{s:Pebbling}

We use the term ``Lagarias-Saks'' pebbling to distinguish it from other forms of graph pebbling (e.g. black pebbling, black-and-white pebbling, etc.) that are used to model and solve other problems in computer science, optimization, and computational geometry.

% =================================
\subsection{Graph theory definitions}
\label{ss:Graphs}

Here we are given a graph $G$ with vertices $V(G)$, edges $E(G)$, and a {\it cost function} $\sw:E(G)\rar\zZ^+$ (the positive integers).
If no edge costs are specified, we default to constant cost 2 everywhere.
We denote the path $P_n=v_1v_2\cdots v_n$ to have $n$ vertices $v_i$ ($i\le n$) and $n-1$ edges $v_iv_{i+1}$ ($i<n$).
In general, we reserve the letter $n=n(G)$ to denote the number of vertices of a graph $G$.

The {\it Cartesian product} $G\Box H$ of two graphs $G$ and $H$ has vertex set $V(G\Box H) = V(G)\times V(H)$ and edge set $E(G\Box H) = \{(u,v_1)(u,v_2)\mid u\in V(G), v_1v_2\in E(H)\}\cup \{(u_1,v)(u_2,v)\mid u_1u_2\in E(G), v\in V(H)\}$.
For example, $P_2\Box P_2$ is isomorphic to the $4$-cycle $C_4$.
We may also recursively define the sequence of $d$-{\it cubes} $\cQ=(Q^1,Q^2,\ldots,Q^d,\ldots)$ by $Q^1=P_2$ and $Q^d = Q^1\Box Q^{d-1} = \Box_{i=1}^d P_2$.
More general $d$-{\it dimensional grids} are defined similarly as a product of $d$ paths: $\Box_{i=1}^d P_{n_i}$, and for any graph $G$ we can define $G^d=\Box_{i=1}^dG$.
The sequence of graphs of interest to us in Section \ref{s:ThreshPebb} is $\cP^d = (P_1^d,P_2^d,\ldots,P_k^d,\ldots)$.

% =================================
\subsection{Pebbling basics}
\label{ss:Pebb}

We are also given a {\it configuration} of pebbles, modeled by a function $C:V(G)\rar\zN$ (the non-negative integers), by which $C(v)$ indicates the number of pebbles on vertex $v$.
The {\it size} of $C$ is defined to be $|C|=\sum_v C(v)$; i.e. the total number of pebbles on $G$.
Additionally, the edges of $G$ are {\it weighted} by a cost function $\sw: E(G)\rar\zN^+$.
Finally, we are given a {\it target} vertex $r$ to which we are challenged to place a pebble, starting from $C$, via a sequence of pebbling steps, which we now describe.

Suppose that $e=uv\in E(G)$ and $C(u)\ge \sw(e)$.
Then the {\it pebbling step} $u\mapsto v$ removes $\sw(e)$ pebbles from $u$ and adds one pebble to $v$.
Thus, after such a step, the resulting configuration $C'$ has $C'(u)=C(u)-\sw(e)$, $C'(v)=C(v)+1$, and $C'(x)=C(x)$ otherwise.
We say that $C$ is $r$-{\it solvable} if it is possible to win this challenge, and $r$-{unsolvable} otherwise.
For example, suppose $r=v_1$ in $P_n$, set $\sw(v_iv_{i+1})=w_i$ for each $1\le i<n$, and let $t=\prod_{i=1}^{n-1}w_i$.
It is straightforward to show by induction that (a) the configuration with $t-1$ pebbles on $v_n$ and 0 elsewhere is $r$-unsolvable, while (b) any configuration of size $t$ is $r$-solvable.

Thus we are led to define the {\it rooted pebbling number} $\p(G,r)$ of a weighted graph $G$ to be the minimum number of pebbles $t$ such that every size $t$ configuration is $r$-solvable.
From the above we have $\p(P_n,v_1)=\prod_{i=1}^{n-1}w_i$.
Next we define the {\it pebbling number} of $G$ as $\p(G)=\max_r\p(G,r)$.
Any configuration of this size is therefore $r$-solvable for every possible target $r$.
Again, it is fairly evident that $\p(P_n)=\p(P_n,v_1)$.

% =================================
\subsection{How pebbling finds zero-sums}
\label{ss:Pebb0Sums}

Since one can read the proof of Theorem \ref{t:LK} in \cite{Chung,ElleHurl}, we only sketch the main idea and connection through an example.

Suppose that we are given $45$ integers, including $32$, $-11$, $31$, $51$, $42$, $-24$, $48$, $75$ and $-15$.
We envision the lattice $L$ of divisors of $45$ (having relation $a\prec b$ when $a|b$), but think of $L$ as a graph $G$. In particular, $G$ has vertices $v_a$ for each divisor $a$ of $45$, with edges between $v_a$ and $v_b$ if (1) $a|b$ and (2) $c\in\{a,b\}$ whenever $a|c$ and $c|b$.
Notice that $G$ is isomorphic to the $2\times 1$ grid $P_3\Box P_2$ because $45=3^2\cdot 5^1$.
The edges of $G$ are then weighted so that the edge $e=v_av_b$, where $a|b$, has cost $\sw(e)=b/a$.

Numbers like $32$, $-11$, and $31$ will be placed as labeled pebbles at vertex $v_1$ because the are relatively prime to $45$; i.e. their $\gcd$ with $45$ equals $1$.
More generally, each number $m$ is placed as a labeled pebbled at vertex $v_k$, where $k=\gcd(m,45)$.
What this placement guarantees is that each pebble by itself satisfies local versions of both properties of interest.
That is, with respect to vertex $v_{15}$, for example, the pebble labeled $30$ is $0\pmod{15}$ and $1/|30|\le 1/|15|$.

Given any three pebbles at $v_1$, such as $32$, $-11$, and $31$, Theorem \ref{t:Geroldinger} guarantees that we can find a subset of them, namely $32$ and $-11$, that is a zero-sum modulo $3$.
Thus we remove them all from $v_1$, create a new pebble labeled $\{32,-11\}$, and place the new pebble at $v_3$.
This represents a pebbling step from $v_1$ to $v_3$ because $\sw(v_1v_3)=3$.
Also, this new pebble satisfies the local properties at $v_3$ because $\gcd(32-11,45)=3$ and $1/|32|+1/|-11|\le 3(1/|1|)= 1/|3|$.
From the summation (not the orders) perspective, the pebble $\{32,-11\}$ acts like a pebble with label $21$ --- we think of $21$ as the {\it nickname} of $\{32,-11\}$.

Now, since the pebbles $51$, $42$, $-24$, and $48$ are also at $v_3$, we can use Theorem \ref{t:Geroldinger} on these four pebbles and the newly arrived ``$21$'' to find a subset, such as $51$, $48$, and $21$, that sums to zero modulo $5$ (which results in zero modulo $15$ because they are each zero modulo $3$ already).
Thus we remove them all from $v_3$, create a new pebble nestedly labeled $\{51,48,\{32,-11\}\}$, and place the new pebble at $v_{15}$.
This represents a pebbling step from $v_3$ to $v_{15}$ because $\sw(v_3v_{15})=5$.
Also, this new pebble satisfies the local properties at $v_{15}$ because $\gcd(51+48+21,45)=15$ and $1/|51|+1/|48|+1/|21|\le 5(1/|3|) 1/|15|$.
The new pebble has nickname $51+48+21=120$.

Finally, since the pebbles $75$ and $-15$ are also at $v_{15}$, we can use Theorem \ref{t:Geroldinger} on these two pebbles and the newly arrived ``$120$'' to find a subset, such as $75$, $-15$, and $120$, that sums to zero modulo $3$ (which results in zero modulo $45$ because they are each zero modulo $15$ already).
Thus we remove them all from $v_{15}$, create a new pebble nestedly labeled $\{75,-15,\{51,48,\{32,-11\}\}\}$, and place the new pebble at $v_{45}$.
This represents a pebbling step from $v_{15}$ to $v_{45}$ because $\sw(v_{15}v_{45})=3$.
Also, this new pebble satisfies the local properties at $v_{45}$ because $\gcd(75-15+120,45)=45$ and $1/|51|+1/|48|+1/|120|\le 3(1/|45|)= 1/|15|$.

The realization that the local properties at $v_{45}$ are equivalent to the sought-after original properties yields the solution $\{75,-15,51,48,32,-11\}$.

The above gives a sense of how Lagarias-Saks pebbling models the sequential construction of a zero-sum sequence with small cross number.
Chung then used retracts to reduce the pebbling problem on divisor lattices (i.e. products of paths) to a similar pebbling problem on cubes (i.e. products of edges), subsequently finding the appropriate pebbling number of cubes.
The consequence of these results is what proves Theorem \ref{t:LK} and, essentially, Theorem \ref{t:Geroldinger}.
There are extra wrinkles to the general abelian group case, involving a specialized representation of them and using Ferrer's diagrams of partitions and their duals.
Furthermore, the technique actually results in the following generalization.
For an additive group $\G$ with subgroup $\sH$, the sequence $g_1, \ldots, g_t$ of elements of $\G$ is an {\it $\sH$-sum} sequence if $g_1 + \cdots + g_t \in \sH$.

\begin{thm}[Elledge-Hurlbert \cite{ElleHurl}]
\label{t:ElleHurl}
For any finite abelian group $\G$ and subgroup $\sH$ of $\G$, every sequence of $|\G|/|\sH|$ of its elements contains an $\sH$-sum subsequence with cross number at most $1/|\sH|$.
\end{thm}

It should be noted that, while Chung's method takes advantage of writing a cyclic group as $\G=\prod_{i=1}^d \zZ_{p_i^{k_i}}$, where $p_1,\ldots,p_d$ are distinct primes, Theorem \ref{t:ElleHurl} allows for any number of these primes to be identical.
In this more general case the graph $G(\G)$ corresponding to the lattice $L=L(\G)$ is not always the divisor lattice of $\prod_{i=1}^d p_i^{k_i}$; instead it is always the product of paths: $G(\G)=\Box_{i=1}^d P_{k_i+1}$.

%%%%%%%%%%%%%%%%%%%%%%%%%%%%%%%%%%%%%%%%%%%%%%%
%       THRESHOLD PEBBLING
%%%%%%%%%%%%%%%%%%%%%%%%%%%%%%%%%%%%%%%%%%%%%%%
\section{Threshold Pebbling} 
\label{s:ThreshPebb}
In this section, we formalize the probabilistic perspective promised in the introduction.
Rather than placing our pebbles on the vertices of a graph carefully, we close our eyes and choose a configuration at random.
We ask ourselves a simple question: ``Is it likely that our random configuration is solvable?"  How many pebbles must we place randomly before we can sleep soundly, confident that our resulting configuration will (probably) be solvable?\\

Of course, this is the same sort of randomness as in the zero-sum sequence framework discussed in the introduction.
We fix some finite abelian group $\Gamma$, and select a random sequence from the group.
How large must this random sequence be, in order to make a zero-sum subsequence likely?  However, our results will follow from threshold pebbling results, so we focus on this framework.
In both this paragraph and the preceding, one may be worried by the use of the imprecise word `likely'.
We now give a formal framework, and will give a precise mathematical meaning to `probably'.

\subsection{Threshold definitions}

We fix a number of pebbles (or sequence length) $t$, and assume that our host graph $G$ (or group $\Gamma$) has order $n$.
We select an initial configuration of $t$ pebbles (sequence elements) uniformly at random among all size $t$ configurations; there are $\binom{n+t-1}{t}$ such configurations.
We denote this probability space by $\mathcal{D}_{G,t}$ ($\mathcal{D}_{\G,t}$).
It is worth noting that this is {\emph{not}} the same distribution that one gets by placing each of $t$ pebbles onto the graph at random -- in this latter scheme, large piles of pebbles are much more likely.
Throughout, we remain interested in the configuration model.
We further note that this is in line with the random model for sequences discussed earlier, where we selected a length $t$ sequence uniformly among all length $t$ sequences in $\G$.

With our probability space in place, we can now begin to make precise our leading question ``How many pebbles are necessary to make it likely that our random initial configuration is solvable?"  In particular, we define the quantity $\Phalf(G)$\footnote{It may seem like our choice of $\nicefrac{1}{2}$ is arbitrary; however, as we'll be focused on cases where our probabilities are either very near 0 or very near 1, this choice will make no difference to the remainder of the work.
Any choice of $p\in(0,1)$ would yield identical theorems.}. Again, this is precisely in line with to our earlier $\Phalf(\cdot)$ definitions in the introduction.
\[\Phalf(G):=\min\left\{t:\pr\left[D\in\mathcal{D}_{G,t}\textrm{ is solvable}\right]\ge\nicefrac{1}{2}\right\}.\]

While determining $\Phalf(G)$ for any particular graph is an interesting problem, it is also a very hard one.
Here, we focus here on sequences of graphs.
That is, given some infinite sequence of graphs $\cG=\left(G_i\right)_{i\in\zN}$, what can we say about the asymptotics of $\Phalf(G_n)$ as a function of $n$?  We again use the the language of thresholds to describe our results.

As an example, let $K_n$ denote the {\it complete graph} on $n$ vertices: every pair of vertices is adjacent.
When our sequence of graphs is the sequence of complete graphs $\cK=\left(K_n\right)_{n\in\zN}$ with uniform edge weights $\sw(e)=2$, then we find ourselves in a situation similar to Feller's Birthday problem (page 33 of \cite{Feller}) -- our configuration is solvable if either we have a pebble on every vertex or if there are two pebbles on any one vertex (since one pebble's sacrifice will allow the survivor to travel to any vertex at all).
The unsolvable $t$-configurations are thus the ones with $t$ pebbles distributed injectively to the $n$ vertices; there are $\binom{n}{t}$ such $t$-configurations, and so the probability that a random configuration is solvable is $1-\binom{n}{t}\big/\binom{n+t-1}{t}$.
%The interested reader should 
It is straightforward to show that this probability tends to $1$ when $t(n)=c(n)\sqrt{n}$ with $c(n)\to\infty$, and that it approaches $0$ if $c(n)\to 0$; this shows that $\sqrt{n}$ is a threshold for $\cK$, though it does not show that this threshold is sharp.
%, but we omit the standard-but-tedious calculations here.
(The only difference between pebbling and birthdays is that in Feller's problem the people are labeled, whereas here the pebbles are unlabeled.)

\subsection{Threshold results}

In \cite{BeBrCzHu}, the authors prove the existence of weak thresholds for monotone families of multisets.
Because the family of unsolvable configurations is monotone (closed under removing pebbles), they obtain the following theorem.

\begin{thm}[Bekmetjev-Brightwell-Czygrinow-Hurlbert \cite{BeBrCzHu}]
\label{t:Exist}
Every infinite sequence of graphs has a weak pebbling threshold function.
\end{thm}

As a concrete example, the sequence of papers \cite{CzEaHuKa,BeBrCzHu,WiSaJaGo} produced ever-sharpening estimates on the ($\sw=2$) pebbling threshold for the sequence of paths, leading to the following result.

\begin{thm}[Czygrinow-Hurlbert \cite{CzygHurlSpec}]
\label{t:CzygHurl}
Let $\cP=(P_1,\ldots,P_n,\ldots)$.
Then for every $c>1$ we have $\t(\cP)\subseteq \W\left(n2^{\sqrt{\lg n}/c}\right)\cap O\left(n2^{c\sqrt{\lg n}}\right)$.
\end{thm}

Note that, since the constant $c$ in the exponent is a multiplicative factor, Theorem \ref{t:CzygHurl} does not yield a weak threshold (although Theorem \ref{t:CzygHurl} guarantees the existence of such a threshold).
Such a threshold was established in a sharp form in \cite{BushKett}.

\begin{thm}[Bushaw-Kettle \cite{BushKett}]
\label{t:BushKettPath} 
With $\cP$ as above, we have 
$\t(\cP)\subseteq\Th\left(ne^{\sqrt{\lg w\lg n}-(\lg\lg n)/2+o(1)}\right)$.
\end{thm}

One can view an $n$-dimensional grid as a product of paths, and so the above results for paths can be seen as the first step toward establishing thresholds for grids of fixed dimension.
A weak threshold was established in \cite{BushKett}, in the somewhat more general setting where pebbling moves in different grid directions can have different costs (we note that the $O(1)$ term in the exponent here as opposed to the $o(1)$ term in the previous theorem is the subtle difference causing one to be weak and one to be sharp.

\begin{thm}[Bushaw-Kettle \cite{BushKett}]
\label{t:BushKettGrid}
Let $P_k^d=\Box_{i=1}^d P_k$, $n=k^d$, and consider the random pebbling model with cost $w_i$ in coordinate $i$.
For fixed dimension $d$, define the sequence $\cP^d = (P_1^d,\ldots,P_k^d,\ldots)$.
Then,
$$\t(\cP^d)= n\exp\left[\left(\frac{(d+1)!\prod_i\lg w_i}{2}\lg k\right)^{\frac{1}{d+1}} - \left(\frac{d}{d+1}\right)\lg\lg k + O(1)\right].$$
\end{thm}

At another extreme, if instead one takes a product of `$P_2$'s with increasing dimension we find the sequence of $d$-cubes, $\cQ$.
Several groups studied thresholds for pebbling $d$-cubes, culminating in the following theorem.

\begin{thm}[Czygrinow-Wagner \cite{CzygWagn}, Alon \cite{AlonThresh}]
\label{t:CzygWagn}
Let $\cQ=(Q^1,\ldots,Q^d,\ldots)$ and $n=2^d$.
Then for every $\e>0$ we have $\t(\cQ)\ \sse\ \W(n^{1-\e})\cap O\left(n/(\lg\lg n)^{1-\e}\right)$.
\end{thm}

In this case, Theorem \ref{t:CzygWagn} translates into the following bounds on the threshold for $\zZ_2^d$.

\begin{thm}
\label{t:Z2d}
Let $\mathcal{Z}=(\zZ_2^1,\zZ_2^2,\ldots)$. Then for every 

Let $Z_d$ be the statement that the sequence $g_1, \ldots, g_t$ in $\zZ_2^d$ is good, and let $\cZ=(Z_1, \ldots, Z_d, \ldots)$ and $n=2^d$.
Then for every $\e>0$ we have
\[\t(\cZ)\ \sse\ O\left(n/(\lg\lg n)^{1-\e}\right).\]
\end{thm}

We end this section with a concrete example of Theorem \ref{t:BushHurl}; we believe that this gives an instructive look at the result in the case of a relatively small group of simple form.

\begin{cor}
\label{t:Example35}
Let $Z_k$ be the statement that the sequence $g_1, \ldots, g_t$ in $\zZ_{35^k}$ is good, and let $\cZ=(Z_1,\ldots,Z_k,\ldots)$.
Then
\[\t(\cZ)\leq n\exp\left[\left(3\lg 5\lg 7\lg k\right)^{1/3} - \left(\frac{2}{3}\right)\lg\lg k + O(1)\right].\]
\end{cor}

%%%%%%%%%%%%%%%%%%%%%%%%%%%%%%%%%%%%%%%%%%%%%%%
%       PROOF OF MAIN THEOREM
%%%%%%%%%%%%%%%%%%%%%%%%%%%%%%%%%%%%%%%%%%%%%%%

\section{Proof of Main Theorem}
\label{s:Proof}
In this section, we use the tools introduced in Sections \ref{s:Pebbling} and \ref{s:ThreshPebb} in order to prove Theorem \ref{t:BushHurl}.  
With this in mind, let $\G\cong \prod_{i=1}^d \zZ_{p_i^k}$ be a finite abelian group, where $p_1,\ldots,p_d$ are (not necessarily distinct) primes and $k\ge 1$.  
Further, let $g=(g_1,\ldots,g_t)$ be an arbitrary length $t$ sequence in $\G$ and define \[F(\G)=n\exp\left[\left(\frac{(d+1)!\prod_i\lg p_i}{2}\lg k\right)^{\frac{1}{d+1}} - \left(\frac{d}{d+1}\right)\lg\lg k + O(1)\right].\]

As described in Section \ref{ss:Pebb0Sums}, the sequence $g$ corresponds to the placement $C$ of $t$ pebbles in the lattice graph $G(\G)$.
By Theorem \ref{t:BushKettGrid} the threshold $\t(\cP^d)$ for the solvability of $C$ equals $F(\G)$.
Hence we know that if $t\gg F(\G)$ then $C$ is solvable with probability tending to 1, and by the discussion of Section \ref{ss:Pebb0Sums} we have that $g$ is good with probability tending to 1. 
\hfill$\Box$

%%%%%%%%%%%%%%%%%%%%%%%%%%%%%%%%%%%%%%%%%%%%%%%
%           COMMENTS
%%%%%%%%%%%%%%%%%%%%%%%%%%%%%%%%%%%%%%%%%%%%%%%
\section{Final Comments}  
\label{s:Comments}

The reason that Theorem \ref{t:BushHurl} is not necessarily an equality is because the converse statement that the pebbling placement $C$ corresponding to a good sequence is not always solvable.
For example, if $g=(1,4)$ in $\zZ_5$ then $g$ is good.
However a pebbling move in $G(\zZ_5)$ requires 5 pebbles, which is impossible, and so $C$ is not solvable.
It still may turn out, however, that the probability that $C$ is solvable, given that $g$ is good, tends to 1, which would yield equality in Theorem \ref{t:BushHurl}.
Even this, though, seems daunting because, for a general placement $C$ of pebbles on an arbitrary graph $G$ with vertex $v$, answering the question ``does $C$ solve $v$?'' is an {\sf NP}-complete problem (see \cite{HurlKier,MilaClar}).
Some respite may be found in the possibility that the question on $d$-dimensional grids is in {\sf P}.

\begin{prob}
\label{p:LowerBounds}
How can one find lower bounds on $\Phalf(\G)$ for any given abelian group $\G$?
\end{prob}

We finish with some remarks and questions regarding other well-known zero-sum problems and conjectures in addition to generalizing Theorem \ref{t:BushHurl}. 

\begin{qst}
\label{q:BushHurl}
Let $\G\cong \prod_{i=1}^d \zZ_{p_i^{k_i}}$ be a finite abelian group, where $p_1,\ldots,p_d$ are (not necessarily distinct) primes and each $k_i\ge 1$.
What is $\Phalf(\G)$?
\end{qst}

Answering Question \ref{q:BushHurl} will likely require generalizing Theorem \ref{t:BushKettGrid}.
Additionally, Theorem \ref{t:BushKettGrid} would also need to be generalized to be able to answer the next question. In this case we have each $k_i=1$, which corresponds to strengthening and generalizing Theorem \ref{t:CzygWagn}.

\begin{qst}
\label{q:Examplepp}
Let $(p_1, p_2, \ldots)$ be an infinite sequence of (not necessarily distinct) primes, $\G_d=\prod_{i=1}^d \zZ_{p_i}$, and $\mathcal{G}=(\G_1,\G_2,\ldots)$.  What is $\tau(\mathcal{G})$?
\end{qst}

For a group $\G$ define the {\it Davenport constant} $\sD(\G)$ to be the smallest $t$ such that every sequence of $t$ elements of $\G$ contains a zero-sum subsequence.
For a finite abelian group $\G$ we write $\G\cong \oplus_{i=1}^k \zZ_{n_i}$, with $n_i | n_{i+1}$ for each $1\le i<k$, and define the {\it Davenport function} $\dav(\G)=\left(\sum_{i=1}^k n_i\right)-k+1$.
In the 1960s, Paul Erd\H{o}s conjectured that every finite abelian group satisfies $\sD(\G)=\dav(\G)$.
While Olson proved that Erd\H{o}s's conjecture was true if $\rank(\G)\le 2$ \cite{OlsonRank2} or if $p$ is prime and $\G$ is a $p$-group $\oplus_{i=1}^k \zZ_{p^{a_i}}$ \cite{OlsonPGroups}, in 1969 van Emde Boas gave the counterexample $\G=\zZ_2^4\oplus\zZ_6$ (for which $\dav(\G)=11>10=\sD(\G)$) and proved subsequently that it is false for ranks four and higher; e.g. $\zZ_3^3\oplus\zZ_6$.
This left the conjecture open for rank 3 groups (now known as Olson's conjecture).
We ask the following questions.

\begin{qst}
\label{q:pebbmodel}
Is it possible to model Olson's conjecture by pebbling (possibly with revised pebbling moves) on some structure?
\end{qst}

Note that there is no cross number condition, and so pebbling on the product of paths is overkill, as evidenced by pebbling numbers that are larger than the Davenport function.

\begin{qst}
\label{q:pebbthresh}
Do (sharp or weak) thresholds exist for zero-sums in this context?
If so, what are these thresholds?
\end{qst}

Let us return to the default pebbling cost $\sw=2$.
For fixed $k$ we define $\cP_k = (P_k^1,P_k^2,\ldots,P_k^d,\ldots)$.
Theorem \ref{t:CzygWagn} shows that $\t(\cP_2)\subseteq o(2^d)$, and it was conjectured in \cite{HurlbertGeneral} that $\t(\cP_k) \subseteq o(k^d)$.
At the other extreme, Theorem \ref{t:BushKettGrid} shows that for fixed $d$ we have $\t(\cP^d)\subseteq \w(k^d)$.
For $d$ a function of $k$ we define $\cP_k^d = (P_1^{d(1)},P_2^{d(2)},\ldots,P_k^{d(k)},\ldots)$.
The following question was asked in \cite{HurlbertGeneral}.

\begin{qst}
[\cite{HurlbertGeneral}]
\label{q:pathprodthresh}
For which function $d=d(k)$ does $\t(\cP_k^d)\subseteq\Th(k^d)$?
\end{qst}

Define the sequence $\cK^d = (K_1^d,K_2^d,\ldots,K_k^d,\ldots)$.
It was shown in \cite{BekmHurl} that $\t(\cK^2) \subseteq \Th(k^{d/2})$.

\begin{qst}
\label{q:cliqueprodthresh}
Is there some $d$ such that $\t(\cK^d) \subseteq \w(k^{d/2})$?
\end{qst}

%%%%%%%%%%%%%%%%%%%%%%%%%%%%%%%%%%%%%%%%%%%%%%%
%           BIBLIO
%%%%%%%%%%%%%%%%%%%%%%%%%%%%%%%%%%%%%%%%%%%%%%%

\bibliographystyle{acm}
\bibliography{refs}

\end{document}